\documentclass[runningheads]{llncs}
\usepackage[T1]{fontenc}
\usepackage{graphicx}
\usepackage{amssymb}
\usepackage{amsmath}
\usepackage{bm}
\usepackage{{booktabs}}
\usepackage{hyperref}

\begin{document}
\title{Learning Deterministic Surrogates for Robust Convex QCQPs}
\author{Egon Peršak\orcidID{0000-0003-3432-2123} \and
        Miguel F. Anjos\orcidID{0000-0002-8258-9116}}
\authorrunning{E. Peršak and M. F. Anjos}
\institute{University of Edinburgh, Edinburgh, UK \\
% \and
% \email{lncs@springer.com}
\email{E.Persak@sms.ed.ac.uk}
}
\maketitle              % typeset the header of the contribution
\begin{abstract}
% The melding of prediction and optimisation pipelines by using differentiable optimisation to optimise for decision losses is a key development for contextual optimisation
Decision-focused learning is a promising development for contextual optimisation. It enables us to train prediction models that reflect the contextual sensitivity structure of the problem. However, there have been limited attempts to extend this paradigm to robust optimisation. We propose a double implicit layer model for training prediction models with respect to robust decision loss in uncertain convex quadratically constrained quadratic programs (QCQP). The first layer solves a deterministic version of the problem, the second layer evaluates the worst case realisation for an uncertainty set centred on the observation given the decisions obtained from the first layer. This enables us to learn model parameterisations that lead to robust decisions while only solving a simpler deterministic problem at test time. Additionally, instead of having to solve a robust counterpart we solve two smaller and potentially easier problems in training. The second layer (worst case problem) can be seen as a regularisation approach for predict-and-optimise by fitting to a neighbourhood of problems instead of just a point observation. We motivate relaxations of the worst-case problem in cases of uncertainty sets that would otherwise lead to trust region problems, and leverage various relaxations to deal with uncertain constraints. Both layers are typically strictly convex in this problem setting and thus have meaningful gradients almost everywhere. We demonstrate an application of this model on simulated experiments. The method is an effective regularisation tool for decision-focused learning for uncertain convex QCQPs.

\end{abstract}
\section{Introduction}
Decision problems are seldom static. Foresighted planning requires us to assert some judgement or prediction of the future to ground the hypothesis space of our decisions. Fortunately, it is common for processes relevant for decision making to have a reasonable degree of predictability either in the form of trends or seasonality. Unfortunately, the dynamics of these processes tend to be stochastic, are often non-linear, potentially non-stationary, and have been observed insufficiently for robust inference. This is a problem as the decision quality is typically quite sensitive to differences between predicted problem parameters and true problem parameters. A small prediction error can render a solution highly suboptimal or even infeasible. 

The classical solution to this from the mathematical optimisation community has been robust optimisation. That is to find a solution that has the best worst case realisation across a neighbourhood of problems defined by an uncertainty set. Determining an appropriate contextual uncertainty set is at least as difficult as the inference problem and is beyond the scope of this work, we assume we have a method which is sufficient for our coverage or regularisation needs. A more recent development that aims to tackle the sensitivity of prediction defined optimisation is decision-focused learning. Instead of attempting to minimise some measure between the predicted and realised problem parameters it aims to minimise the resulting decision loss. The learnt models reflect the underlying contextual sensitivity structure of the problem. The predictions effectively trick the optimisation model into making more robust solutions given the sparsely observed contextual realisations. 

There are two key shortcomings with the decision-focused learning approach. It is computationally demanding as computing the gradients for decision losses requires differentiating through an $\mathrm{arg}\min$ which typically requires at least one computation of an optimal solution. In this work we focus on convex optimisation problems for which we typically rely on differentiating the optimality conditions for gradient computation, but approximate gradient methods have even been developed for blackbox solvers for combinatorial problems \cite{VlastelicaEtal2020:BBoxSolvers}. Decision-focused learning therefore scales poorly with increasing levels of the computational difficulty of optimisation problems. It will scale better for linear programs than quadratic programs, better for quadratic than semidefinite programs (SDP), better for SDP than integer and other non-convex problems. 

The prediction problem for optimisation is typically very high dimensional, as we need to predict the vectors and matrices that define an optimisation problem. As observations are limited, this means we are necessarily dealing with a lot of empirical uncertainty on top of aleatoric uncertainty and are faced with a few-shot learning problem. While standard decision-focused learning ameliorates the sensitivity issue, it does not resolve it without having many observations in similar contexts. Even in data-rich cases the underlying aleatoric uncertainty means that predictions based on minimising empirical regret will overfit to observed realisations for a local context. A potential solution is taking the decision losses from the robust problem centred on the observation. In its classical form this exacerbates the computational issue as the robust counterpart of an optimisation problem is at best (with carefully designed uncertainty sets) in the same class of computational difficulty, but larger. In general, uncertainty sets have to be carefully designed to preserve tractability and robust counterparts usually belong to more computationally challenging problem families \cite{ben1998robust}. The nature of the decision problem may mean that we have insufficient time to obtain a robust solution at test time.

To address these shortcomings we propose a model which learns parameterisations for surrogate deterministic programs whose solutions approximate robust solutions. This means we have to separate the deterministic version of the program from the worst case problem. In the worst case problem we are trying to find a combination of the worst case outcome and constraint violation for a given decision. It acts as a component of the robust decision loss function and thus does not have to be computed at test time. We interpret it as a regularisation method for decision-focused learning. In effect, we are trying to solve the inverse problem of what parameters for a deterministic problem would achieve robust solutions for a given uncertainty set.

We narrow our attention to uncertain convex QCQP. This research focus is motivated by two properties. The first is that for strictly convex QCQPs the output in terms of parameters is continuous and thus subdifferentiable everywhere and differentiable almost everywhere. The second is that QCQPs may have optimal solutions in the interior of the feasible set and therefore the surrogate models can either learn in the direction of the constraints or objective. In contrast, linear programs have outputs which are piecewise constant requiring some form of regularisation in the objective function for decision-focused learning, for example a 2-norm on the decision vector \cite{wilder2019melding} to compute meaningful gradients. This turns the problem into a second order cone program (SOCP), which is the complexity for robust linear programming with ellipsoidal uncertainty. Furthermore, the solution to a robust linear program may not be at a vertex meaning a linear surrogate will potentially have trouble approximating it, or in the case of fixed constraints will not be able to approximate it at all. We are confident our approach can be extended to more general convex optimisation frameworks, but it is beyond the scope of this work. Since general QCQPs tend to be tackled with relaxations that reformulate them into SDP extending robust predict-and-optimise to uncertain SDPs is a promising future direction.

The remainder of this work is structured as follows. Section \ref{background} is a non-exhaustive overview of the prior work in decision-focused learning, robust convex optimisation, and online methods for QCQPs. Section \ref{main} delineates the contribution of this work by presenting the two layer approach for deterministic surrogates and why it works. Section \ref{experiments} presents synthetic experiments demonstrating the usefulness of learnt deterministic surrogates. The work is concluded with a discussion of potential applications, and future directions research in section \ref{conclusion}.

\section{Background}
\label{background}
% draft
% mention robust losses, surco
The field of contextual stochastic optimisation concerns itself with decision making under uncertainty in a data driven way. A paradigm which has seen much excitement in this community is integrating the prediction and optimisation step, which has been variously termed as predict-and-optimise, decision-focused learning, and end-to-end learning. From a differentiable programming standpoint it sees optimisation problems as just another differentiable component in a mapping of data context to decisions. The field of deep learning sees it as an example of an implicit layer. 

Three main directions have been developed to extract gradients from optimisation layers: using the implicit function theorem to differentiate the optimality conditions \cite{donti2017task,amos2017optnet,wilder2019melding,agrawal2019differentiating}, using a mathematically sensible surrogate loss function \cite{elmachtoub2022smart,VlastelicaEtal2020:BBoxSolvers}, or using a surrogate differentiable optimiser \cite{berthet2020learning,kong2022end}. A host of techniques have been developed and they differ in what types of uncertain optimisation problems they are applicable for; for an overview we recommend section 5 of \cite{sadana2023survey}. Notably few deal with constraint uncertainty. Differentiable optimisation has been leveraged to learn linear objective surrogate models for non-linear objective combinatorial optimisation models \cite{ferber2023surco}, which are otherwise largely intractable. The non-linear objective effectively acts as a loss function for the decisions produced by the surrogate model. Robust decision losses were first proposed by \cite{schutte2023robust}. Their analysis shows that for decision-focused learning applied to uncertain optimisation problems empirical regret does not always equal expected regret and is therefore an ineffective loss surrogate. Though the analysis is done for linear programs with uncertain objectives, the same argument can be extended for general uncertain convex optimisation problems. They propose three losses that are better approximates of expected regret. Using a surrogate loss (SPO+ \cite{elmachtoub2022smart}, Perturbed Fenchel–Young Loss \cite{berthet2020learning}) function approach they modify the target decision given the observation $x^*(c)$ to the optimal decision to the robust problem $x_{RO}^*(c)$. This requires a pass through the robust problem in training, but they choose an uncertainty set for which the robust counterpart is at the same level of complexity. Our approach extends the concept of robust decision losses to uncertain convex QCQPs using an implicit differentiation approach.

QCQP are optimisation problems where the objective and all constraint functions $\leq0$ are quadratics of the form $f(\textbf{x})=\textbf{x}^TA\textbf{x} + \textbf{b}^T\textbf{x}+c$. If $A$ is positive semidefinite (PSD) the problem is convex and efficiently solvable using the interior point method for SOCP \cite{nesterov1994interior}. In the special case of $A=0$ this reduces to a linear programming problem. QCQPs have applications in fields such as machine learning, finance, robotic control, and energy systems.

There are well known results for uncertain QCQP using robust optimisation. A robust QCQP with ellipsoidal uncertainty has an exact SDP equivalent \cite{ben1998robust}. There are convex uncertainty sets for which the robust counterpart is known to be NP-hard, such as an intersection of ellipsoids. For large robust QCQP the resulting SDP formulation may be too large for existing solvers to handle so a host of methods have been developed that iteratively solve the deterministic problem to approximately solve the robust one \cite{ben2015oracle,kroer8performance}. Our approach is of similar iterative nature and in the case of only having one context it can be seen as an approximate solution method for robust QCQPs.

\section{Learning Deterministic Surrogates}
\label{main}
\subsection{The Optimisation Problem}
We propose a double implicit layer strategy for learning deterministic surrogates for robust QCQPs. These problems take the form of:
\begin{equation}
\label{reference_problem}
\begin{array}{rlcl}
\displaystyle \textbf{x}^* \in \mathrm{arg}\min_{\textbf{x}} & \textbf{x}^{T}Q\textbf{x} + \textbf{c}^T \textbf{x} + q \\[4pt]

\textrm{s.t.} &  \textbf{x}^{T}A_i\textbf{x} + \textbf{b}_i^T \textbf{x} + \gamma_i \leq 0 & \quad \forall i \in {1,...,n} 
\\[4pt]

& (A_i,\textbf{b}_i,\gamma_i) \in \mathcal{U}_i & \quad \forall i \in {1,...,n} 
\end{array}
\end{equation}
Unlike most predict-and-optimise work where the uncertainty is only in the objective we assume the uncertainty is only in the constraints and that it is constraint-wise. This is more general as a problem with an uncertain objective can trivially be reformulated as an equivalent problem with a certain objective and uncertain constraints. The constraints in this problem can be equivalently written as: 
\begin{equation}
\label{worst_case_constraint}
   \max_{(A_i,\textbf{b}_i,\gamma_i) \in \mathcal{U}_i} \textbf{x}^{T}A_i\textbf{x} + \textbf{b}_i^T \textbf{x} + \gamma_i \leq 0 \quad \forall i \in {1,...,n} 
\end{equation}
We assume $Q,A_i$ are positive semi-definite, though as we later point out this assumption can be relaxed for $A_i$ in our method for certain types of uncertainty sets. We want to determine a set of problem parameters $$\hat{{\mathcal{P}}} = \{\hat{Q},\hat{\textbf{c}},\hat{q},(\hat{A}_i,\hat{\textbf{b}}_i,\hat{\gamma}_i) | \forall i \in {1,...,n} \}$$ such that fixing the parameters in problem (\ref{reference_problem}) yields an optimal solution $\textbf{x}_{det}^*$ which satisfies all worst case constraints (\ref{worst_case_constraint}) and approximately minimises the objective. There purpose of this is two-fold: for a single problem it represents a solution method for robust QCQPs with convex uncertainty sets, and more importantly, when learned with context $\textbf{z}_j$ it enables decision-focused learning to predict problem parameterisations robustly.

Given some context $\textbf{z}_j$ we want to predict a set of parameters $\hat{{\mathcal{P}}}_j$ that yields a decision $\textbf{x}_{det,j}^*$ which is feasible and minimises a version of problem (\ref{reference_problem}). Note that we do not predict $\hat{Q}$ or $\hat{A}_i$ directly, instead we predict some matrix of the same size $\hat{M}_i$ and set the problem values as $\hat{M}_i^T\hat{M}_i$ to ensure they are PSD. The version of the problem is based on observed realisations of parameters ${\mathcal{P}}_j$ and some fixed or conditional uncertainty set ${\bm{\mathcal{U}}_j({\mathcal{P}}_j)}$ centred on the realisation. How the uncertainty set is determined is beyond the scope of this work, but we assume that a robust method is available. The shape and scale of the uncertainty set can be understood as determining the degree of regularisation of the uncertain optimisation problem.

\subsection{Double Implicit Layer Model}
\subsubsection{Predict-and-Optimise} For this prediction task we propose a two block neural network (NN) composed of a prediction block (PR) parameterised by parameters $\theta_P$ and a robust optimisation block (ROB). The prediction block can consist of any standard NN components, so long as it is capable of multioutput prediction (since we have to predict vectors and matrices in this case).
\begin{equation}
\mathrm{PR}(\textbf{z}_j;\theta_P) = \hat{\mathcal{P}}_j
\end{equation}
ROB consists of two implicit layers. The first layer $\mathrm{ROB}_1$ solves a deterministic version of the problem where the inputs are the fixed/predicted parameters from the prediction block and outputs the corresponding optimal decisions.
\begin{equation}
\begin{array}{rlcl}
\displaystyle \mathrm{ROB}_1(\hat{\mathcal{P}}_j): \quad
\textbf{x}_{det,j}^* \in \mathrm{arg}\min_{\textbf{x}} & f(x,\hat{{\mathcal{P}}}_j) 
\\[4pt]
\textrm{s.t.} &  g_i(x,\hat{{\mathcal{P}}}_j) \leq 0 & \quad \forall i \in {1,...,n} 
\end{array}
\end{equation}
where $f,g_i$ are the corresponding quadratic objective and constraint functions written concisely for convenience. We use the belongs to notation since the optimal solution may not be unique, but in practice the layer will output a unique optimal solution. In the standard predict-and-optimise set-up this output would be used to compute a variant of decision loss, typically some regret based metric like SPO: $f(\mathrm{ROB}_1(\hat{{\mathcal{P}}}_j),\mathcal{P}_j) - f(\mathrm{ROB}_1({\mathcal{P}}_j),\mathcal{P}_j)$. 

\subsubsection{Worst Case Layer:}
In ROB we augment this by adding a second implicit layer $\mathrm{ROB}_2$ which evaluates the worst case outcome given an uncertainty set. Since the uncertainty is in the constraints we want to penalise constraint violation in the objective function. Designing constraint penalties depends on the form of the uncertainty set. The worst case problem for robust QCQPs is known to be tractable for simple ellipsoid uncertainty sets. We look at two forms of defining ellipsoid uncertainty sets for the PSD quadratic term coefficient $A_i = P_i^TP_i$. Note that PSD matrices can always be expressed in this form using Cholesky factorisation. The uncertainty sets can either be defined in terms of $A_i$ or $P_i$.

The quadratic functions that describe the constraints of the problem are of the form $\textbf{x}^{T}A_i\textbf{x} + \textbf{v}^T \textbf{x} + s$. The first term can be restated in a lifted form $A_i \bullet \textbf{X}$ where $\textbf{X}=\textbf{x}\textbf{x}^T$, and $\bullet$ represents the Frobenius inner product $\mathrm{tr}(A_i^T \textbf{X})$, the sum of the element-wise product. This means that if the uncertainty set is convex and defined in terms of $A_i$, for fixed $\textbf{x}$ the objective function is linear. However, a problem arises in the case of defining ellipsoidal uncertainty sets in terms of $P_i \in \{ P_i | P_i = P_{i,0} + \sum_{i=k}^{n}u_kP_{i,k},||\textbf{u}||_2 \leq 1 \}$. The worst case constraint (\ref{worst_case_constraint}) becomes a maximisation of a convex function problem, also known as a trust region problem, which is non-convex and therefore generally intractable. However, this problem has some hidden convexity and can be restated as an equivalent SDP using the S-lemma \cite{ben1998robust,bertsimas2011theory}. The worst case constraint in (\ref{worst_case_constraint}) is equivalent to the existence of a scalar $l$ such that a specific matrix $S_i$ composed of affine transformations of the elements defining the uncertainty set, $l$, and the decision vector is positive semi-definite. 

\subsubsection{Case $P_i$:} Our formulation of $\mathrm{ROB}_2$ for the case of $P_i$ converts the linear matrix inequality into a concave penalty. We fix the decision vector to the output of $\mathrm{ROB}_1$, so the only variable to determine the matrix $S_i$ is $l$. For the decision to be feasible we need to determine if an $l$ exists such that $S_i$ is PSD. This is equivalent to showing that the smallest eigenvector $\lambda_{min}$ of $S_i$ is non-negative. The function $\lambda_{min}$ is concave and tractable in a maximisation problem. Denote the collection of $S_i$ matrices for each of the constraints as $\textbf{S}$ and the corresponding vector of $l_i$ as $\textbf{l}$ where a pair $S_i,l_i$ reflects the feasibility condition components for each constraint i. We propose the following problem to determine the degree of worst case constraint violation:
\begin{equation}
\label{rob2_pi}
\begin{array}{rlcl}
\displaystyle
\mathrm{ROB}_2(\textbf{x}_{det,j}^*,{\mathcal{P}}_j;\bm{\mathcal{U}}_{j}) : \textbf{S}_j^* \in \mathrm{arg}\max_{\textbf{l},\textbf{S}} 
& \sum_{i=1}^n\lambda_{min}(S_i)
\\[4pt]
\textrm{s.t.} & S_i = S({\bm{\mathcal{U}}_{j,i}({\mathcal{P}}_j),\textbf{x}_{det,j}^*,{l}_i})\\[4pt]
& {l}_i \in \mathbb{R}\\[4pt]
& \forall i \in {1,...,n} 
\end{array}
\end{equation}
Where $S(\cdot)$ is the linear matrix function that composes the parameters from the uncertainty set, the fixed decisions, and the $l$ variable into the feasibility condition matrix. If the problem is feasible for $\textbf{x}_{det,j}^*$ then the smallest eigenvalue of each of the solutions $\textbf{S}_j^*$ in (\ref{rob2_pi}) is greater than 0.  If the fixed decision is infeasible $\textbf{S}_{j,i}^*$ reflects how close to feasibility the fixed decision is via its eigenvalues. We can therefore use it as a penalty for uncertain constraint violation. The fixed decision is also feasible if (\ref{rob2_pi}) is unbounded and we manage this exception by setting all penalties to 0. The loss function in this case is:
\begin{equation}
\label{loss_1}
    \mathcal{L}(\textbf{x}_{det,j}^*,\textbf{S}_j^*) = f(\textbf{x}_{det,j}^*,\mathcal{P}_{j}) + \sum_{i=1}^n \left[ \tau_i (\sum\mathrm{ReLU}(
    -\lambda(\textbf{S}_{j,i}^*))) \right]
\end{equation}
Where $\lambda(\cdot)$ is a differentiable function $\lambda : \mathbb{R}^{m \times m} \rightarrow \mathbb{R}^{m}$ that returns the eigenvalues of of a matrix and $\tau_i \geq 0$ is a penalty coefficient. The right term penalises any negative eigenvalues in $\textbf{S}_j^*$ which reflects constraint violation. If there is no constraint violation gradients will only be taken with respect to the deterministic problem. In both cases at test time we simply deactivate $\mathrm{ROB}_2$ to get the deterministic surrogate solution.

\subsubsection{Case $A_i$:} Our formulation of $\mathrm{ROB}_2$ for the case of an $A_i$ ellipsoid is simpler. For a fixed $\textbf{x}_{det,j}^*$ the worst case problem has a linear objective function with a SOCP constraint. The robust counterpart for such a problem is also easier (SOCP) than in case $P_i$ (SDP). We reflect the magnitude of constraint violation by applying Lagrangian relaxation to the uncertain constraints. We assign a penalty coefficient of $\lambda_i < 0$ to the value of each constraint as we want to maximise the violation. The robustness penalty problem for this case is:
\begin{equation}
\label{rob2_ai}
\begin{array}{rlcl}
\displaystyle
\mathrm{ROB}_2(\textbf{x}_{det,j}^*,{\mathcal{P}}_j;\bm{\mathcal{U}}_{j}): \mathcal{P}_{wc}^* \in \mathrm{arg}\max_{\mathcal{P}_{wc}} &  \sum_{i=1}^n\lambda_i g_i(\textbf{x}_{det,j}^*,\mathcal{P}_{wc,i})
\\[4pt]
\textrm{s.t.} & \mathcal{P}_{wc,i} \in  {\bm{\mathcal{U}}_{j,i}({\mathcal{P}}_j)}  \quad \forall i \in {1,...,n}
\end{array}
\end{equation}
The decision variables in this case are the uncertain matrix, vector, and scalar parameters that define the QCQP. The optimal solution to this problem is the worst case realisation in terms of constraint violation. The corresponding loss function is:
\begin{equation}
\label{loss_2}
    \mathcal{L}(\textbf{x}_{det,j}^*,\mathcal{P}_{wc,j}^*) = f(\textbf{x}_{det,j}^*,\mathcal{P}_{j}) + \sum_{i=1}^n\mathrm{ReLU}(\lambda_i  g_i(\textbf{x}_{det,j}^*,\mathcal{P}_{wc,j}^*))
\end{equation}
The penalty for constraint violation in the objective has one obvious flaw, it awards being more in the interior of the feasible set. To correct for that we apply a ReLU activation function to the penalty terms so that the constraint violation penalty is only applied if the current solution is infeasible for that constraint.

The simplicity of case $A_i$ warrants a rethinking of how uncertainty sets $\bm{\mathcal{U}}$ should be defined for this method for larger instances. We can define the uncertainty set as any convex set that is compatible with the subfamilies of conic programming. In this setting, given that we have a linear objective we do not necessarily have to restrict ourselves to PSD quadratic coefficient terms. This means that we can learn convex QCQP surrogates for uncertain problems that have potentially non-convex quadratic realisations. Moreover, defining uncertainty sets within this framework does not require specialised mathematical understanding of how to tractably reformulate a problem.

\subsection{The Learning Problem}
The backwards pass relies on finding the gradient of the loss function with respect to $\theta_P$. Using the chain rule we can decompose gradient calculation as:
\begin{equation}
    \frac{\partial\mathcal{L}}{\partial \theta_P} = 
     \left[ \frac{\partial \mathcal{L}}{\partial \textbf{x}_{det,j}^*}
     + 
    \frac{\partial \mathcal{L}}{\partial \mathrm{ROB}_2^*}
    \frac{\partial \mathrm{ROB}_2(\textbf{x}_{det,j}^*,{\mathcal{P}}_j;\bm{\mathcal{U}}_{j,i})}{\partial \textbf{x}_{det,j}^*}
    \right]
    \frac{\partial \mathrm{ROB}_1(\hat{{\mathcal{P}}}_j)}{\partial \hat{{\mathcal{P}}}_j}
    \frac{\partial\mathrm{PR}(\textbf{z}_j;\theta_P)}{\partial \theta_P}
\end{equation}
There are two terms in this expression that are not easy to evaluate: the partial derivatives for the $\mathrm{ROB}$ layers. These can be obtained using an appropriate method for differentiable optimisation. In our case, since we are dealing with conic representable problems with uncertain constraints the method of differentiating through a conic program \cite{agrawal2019differentiating} is appropriate. It obtains the derivative by implicitly differentiating the conic programme's homogeneous self-dual embedding and is conveniently implemented in the package cvxpylayers \cite{cvxpylayers2019}. Note that cvxpylayers has been observed to struggle with large problems as differentiating the KKT conditions has worst-case cubic complexity \cite{sun2022alternating} in the number of variables and constraints. ADMM has been adapted to efficiently deal with large differentiable quadratic programs \cite{butler2023efficient}. We refer the reader to \cite{sadana2023survey} for an overview alternatives.

While a differentiation method for conic representable optimisation problems exists, it does not mean that we have meaningful gradients everywhere. In the case of linear programming the optimal value as a function of problem parameters is a piece-wise constant discontinuous function. This means that we have a gradient of 0 at any point that is differentiable. We know that for strictly convex QP the output in terms of parameters is continuous and theorem 1 from \cite{amos2017optnet} shows that this means that the implicit function is differentiable almost everywhere and subdifferentiable everywhere. Quadratic programs can be equivalently reformulated as linear objective QCQPs. In our setup $\mathrm{ROB}_1$ is a QCQP, assuming that all coefficient matrices are PD making the problem strictly convex theorem 1 applies to $\mathrm{ROB}_1$.
% FINISH ARGUMENT
% \begin{theorem}

% \end{theorem}
In $\mathrm{ROB}_2$ the nature of the implicit function depends on the definition of the uncertainty set. In the case $A_i$ all constraint functions of the problem are strictly convex with a linear objective. This means that a change in the direction or the orientation of the linear objective will result in a continuous change of the optimal solution. In case $P_i$ the optimisation problem is to maximise a smallest eigenvalue function and a set of linear matrix equality constraints describing the matrix, this is equivalent to minimising the lowest eigenvalue of the negative of the matrix, which is known to be a convex problem that can be cast as an SDP \cite{mengi2014numerical}. The optimal output of SDPs is typically continuous in terms of problem parameters and should have meaningful gradients almost everywhere.

\section{Experiments}
\label{experiments}
We conducted two small synthetic experiments to demonstrate a proof of concept\footnote{Preliminary code available at: \\ \href{https://github.com/EgoPer/Deterministic-Surrogates-for-Uncertain-Convex-QCQP}{https://github.com/EgoPer/Deterministic-Surrogates-for-Uncertain-Convex-QCQP}}. The first experiment shows that the algorithm approximates the robust counterpart solution on a single, no-context instance problem with a case $P_i$ uncertainty set. The second is an application to a simulated contextual optimisation problem with a case $A_i$ uncertainty set where we compare the results with a non-regularised predict-and-optimise method.

\subsubsection{Experiment 1} We replicate problem (\ref{reference_problem}) with an ellipsoid uncertainty set for $P_i,\textbf{b}_i,\gamma_i$ where $P_i^TP_i=A_i$. The quadratic coefficient in the objective $Q$ and scalar $q$ were set to zero. This recreates the problem studied in \cite{ben1998robust} for which there is a widely known tractable convex robust counterpart (RC). We randomly generate the problem coefficients at five problem sizes with five constraints, each with an uncertainty set defined by four points. It is first solved using the robust counterpart. We randomly initialise the deterministic surrogate problem (SUR) and apply gradient descent in the form of Adam for two hundred steps. We select the feasible solution with the lowest objective value as the output of this process. During optimisation the solution $\textbf{x}_{det,j}^*$ oscillates between being feasible (no constraint violation penalty gradients) and infeasible. Figure \ref{exp1} shows the narrow gap between the known RC optimal solution and the approximation using robust surrogates.
\begin{figure}
    \centering
    \begin{tabular}{c|c|c|c}
     Problem size &  RC opt value &  SUR opt value &  gap after 200 \\
    \midrule
               10 &               -0.1949 &                  -0.1816 &              0.0683 \\
               20 &               -0.5704 &                  -0.5402 &              0.0529 \\
               30 &               -0.5077 &                  -0.4809 &              0.0528 \\
               40 &               -0.5866 &                  -0.5479 &              0.0660 \\
               50 &               -0.6000 &                  -0.5608 &              0.0653 \\
    \end{tabular}
    \caption{Comparing Deterministic Surrogates with Robust Counterparts}
    \label{exp1}
\end{figure}
So far we have evidenced that our algorithm approximately converges to an exact algorithm two orders of magnitude slower. What is promising is that $\mathrm{ROB}_1$ for this problem is resolved two orders of magnitude faster than the robust counterpart, which reflects the difference in computational difficulty between QCQP and SDP.

\subsubsection{Experiment 2} The real benefit of this method is in contextual settings where one may be inclined to apply decision-focused learning. In experiment 2 we generated a contextual QCQP. We  generate a fixed correlation matrix $C$ based on randomly generated positive eigenvalues. We generate context $z_{j,i} \sim \mathrm{Unif}(0.1,1)$ where $j$ indexes over the hundred samples we generate and $i$ over the number of contextual covariates (in our case four). We then generate conditional variances $\sigma_{j,i}^2$ as convex combinations of the context $z_{j,i}$, the weights for which are also randomly generated. Let $D_j$ be a matrix for which the diagonal are the generated conditional variances. We define conditional covariance matrices as $\Sigma_j = D_j^{-\frac{1}{2}}C D_j^{-\frac{1}{2}}$. These conditional covariance matrices are used to parameterise a conditional Wishart distribution with a degree of freedom of fifty from which we sample the uncertain PSD matrices for our problem and then divide them by the degree of freedom. These values should reflect the distribution of empirical covariance matrix estimates we would observe for an underlying process generated by $\mathcal{N}(0,\Sigma_j)$. Our optimisation problem is:
\begin{equation}
\begin{array}{rlcl}
\displaystyle \quad
\textbf{x}_{det,j}^* \in \mathrm{arg}\min_{\textbf{x}} & -\textbf{c}^T\textbf{x}
\\[4pt]
\textrm{s.t.} & \textbf{x}^T\Sigma_j \textbf{x} \leq r \\[4pt]
& \textbf{x}^T\mathbf{1} \leq 1, \quad \textbf{x} \geq 0
\end{array}
\end{equation}
Where $\textbf{c} \geq 0$ is a fixed cost vector and $r$ is some defined level of acceptable variance-denominated risk. A possible interpretation of this problem is an investment manager deciding what bonds to hold to maturity to maximise yield while making sure risk levels do not force liquidation. Since $\Sigma_j$ is uncertain a deterministic solution to this problem could be infeasible. 

We train a simple two layer neural network to predict $\hat{c}_j,\hat{\Sigma}_j$ based on input $z_j$ using decision losses and robust decision losses with the method of deterministic surrogates. We define our regularising uncertainty set as:
$$\mathcal{U}_j(\Sigma_j) = \{ \Sigma_j + A | A=A^T,||A||_{F} \leq 1\}$$
The two prediction models are trained on seventy samples and tested on thirty. Note that we use loss function (\ref{loss_2}) in both cases, however, the worst case outcome value is replaced with the observed covariance in the vanilla decision loss case. In our initialisation the robust decision loss produces predictions which cause feasible decisions in all thirty test cases, whereas decision loss alone produces predictions which cause decisions which are infeasible in thirteen out of the thirty cases. The average optimal value without regularisation is $37.3\%$ lower, but infeasible for almost half of the cases. These values are obviously very sensitive to initialisation, but in general we observe that robust decision losses yield more feasible test cases. We are excited about further experimentation to exhibit the merits of robust decision losses for QCQP and beyond.

\section{Conclusion}
\label{conclusion}
We establish a method for obtaining robust decision losses for uncertain convex QCQP. By casting the problem as trying to determine a contextual deterministic surrogate we obtain gradients with respect to the worst case in a neighbourhood of problem realisations. 

The principal advantage of this method is the speedup at test time achieved by separating the robust problem into two optimisation problems. At test time it means we can approximate robust solutions in the time it takes to optimise deterministic problems. At train time instead of solving one larger robust counterpart, two smaller ones are solved. $\mathrm{ROB}_1$ is a QCQP, whereas $\mathrm{ROB}_2$ is in the same family as the robust counterpart, but with fewer decision variables which can be faster to solve on average. The method can be seen as regularisation for decision focused learning.

We see promising research directions in exploring how different choices of uncertainty sets for robust losses affect decision-focused learning for QCQPs, in extending the methodology to uncertain SDPs, and experimenting with the method on real-world problems.

\begin{credits}
% \subsubsection{\ackname} A bold run-in heading in small font size at the end of the paper is
% used for general acknowledgments, for example: This study was funded

\subsubsection{\discintname}
The authors have no competing interests to declare that are
relevant to the content of this article.
\end{credits}

\bibliographystyle{splncs04}

\bibliography{bib_paper}

\end{document}